\documentclass[12pt]{article}
\usepackage{mathrsfs}
\usepackage{amscd}
\usepackage{amsmath,amsfonts,amssymb,amscd}
\usepackage{indentfirst,graphics,epsfig,psfrag}
\input{epsf}
\usepackage{ifpdf}
\usepackage{enumerate}
\usepackage{appendix}
\usepackage{enumerate}

\usepackage{lineno}

\makeatletter \@addtoreset{figure}{section} \makeatother
\makeatletter
\long\def\@makecaption#1#2{%
   \vskip 10\p@
   \setbox\@tempboxa\hbox{{#1}\ \ #2}%
   \ifdim \wd\@tempboxa >\hsize
       {#1}\ \ #2\par
   \else
       \hbox to\hsize{\hfil\box\@tempboxa\hfil}%
   \fi}
\makeatother

\newtheorem{thm}{Theorem}[section]
\newtheorem{cor}[thm]{Corollary}
\newtheorem{lem}[thm]{Lemma}

\newtheorem{obs}[thm]{Observation}
\newtheorem{pro}[thm]{Proposition}

\newcommand{\qed}{{\hfill\rule{3pt}{7pt}}}

\def\qed{\hfill \rule{4pt}{7pt}}

\begin{document}
\title{On strong rainbow connection number\footnote{Supported by NSFC.}}
\author{
\small  Xueliang Li, Yuefang Sun\\
\small Center for Combinatorics and LPMC-TJKLC\\
\small Nankai University, Tianjin 300071, P.R. China\\
\small E-mails: lxl@nankai.edu.cn, syf@cfc.nankai.edu.cn
 }
\date{}
\maketitle
\begin{abstract}
A path in an edge-colored graph, where adjacent edges may be colored
the same, is a rainbow path if no two edges of it are colored the
same. For any two vertices $u$ and $v$ of $G$, a rainbow $u-v$
geodesic in $G$ is a rainbow $u-v$ path of length $d(u,v)$, where
$d(u,v)$ is the distance between $u$ and $v$. The graph $G$ is
strongly rainbow connected if there exists a rainbow $u-v$ geodesic
for any two vertices $u$ and $v$ in $G$. The strong rainbow
connection number of $G$, denoted $src(G)$, is the minimum number of
colors that are needed in order to make $G$ strong rainbow
connected. In this paper, we first investigate the graphs with large
strong rainbow connection numbers. Chartrand et al. obtained that
$G$ is a tree if and only if $src(G)=m$, we will show that
$src(G)\neq m-1$, so $G$ is not a tree if and only if $src(G)\leq
m-2$, where $m$ is the number of edge of $G$. Furthermore, we
characterize the graphs $G$ with $src(G)=m-2$. We next give a sharp
upper bound for $src(G)$ according to
the number of edge-disjoint triangles in graph $G$, and give a necessary
and sufficient condition for the equality. \\[2mm]
{\bf Keywords:} edge-colored graph, rainbow path, rainbow geodesic,
strong rainbow connection number, edge-disjoint triangle.\\[2mm]
{\bf AMS Subject Classification 2000:} 05C15, 05C40
\end{abstract}

\section{Introduction}

All graphs in this paper are finite, undirected and simple. Let $G$
be a nontrivial connected graph on which is defined a coloring
$c:E(G)\rightarrow \{1,2,\cdots,n\}$, $n\in \mathbb{N}$, of the
edges of $G$, where adjacent edges may be colored the same. A path
is a $rainbow$ $path$ if no two edges of it are colored the same. An
edge-coloring graph $G$ is $rainbow ~connected$ if any two vertices
are connected by a rainbow path. Clearly, if a graph is rainbow
connected, it must be connected. Conversely, any connected graph has
a trivial edge-coloring that makes it rainbow connected; just color
each edge with a distinct color. Thus, we define the
$rainbow~connection~number$ of a connected graph $G$, denoted
$rc(G)$, as the smallest number of colors that are needed in order
to make $G$ rainbow connected. Let $c$ be a rainbow coloring of a
connected graph $G$. For any two vertices $u$ and $v$ of $G$, a
$rainbow$~$u-v$~$geodesic$ in $G$ is a rainbow $u-v$ path of length
$d(u,v)$, where $d(u,v)$ is the distance between $u$ and $v$. The
graph $G$ is $strongly$~$rainbow$~$connected$ if there exists a
rainbow $u-v$ geodesic for any two vertices $u$ and $v$ in $G$. In
this case, the coloring $c$ is called a
$strong$~$rainbow$~$coloring$ of $G$. Similarly, we define the
$strong~rainbow~connection~number$ of a connected graph $G$, denoted
$src(G)$, as the smallest number of colors that are needed in order
to make $G$ strong rainbow connected. A strong rainbow coloring of
$G$ using $src(G)$ colors is called a
$minimum~strong~rainbow~coloring$ of $G$. Clearly, we have
$diam(G)\leq rc(G)\leq src(G)\leq m$ where $diam(G)$ denotes the
diameter of $G$ and $m$ is the size of $G$. In an edge-colored graph
$G$, we use $c(e)$ to denote the color of edge $e$, then for a
subgraph $G_1$ of $G$, $c(G_1)$ denotes the set of colors of edges
in $G_1$.

In \cite{Y. Caro}, the authors investigated the graphs with small
rainbow connection numbers, they determined a sufficient condition
that guarantee $rc(G)=2$.

\begin{thm}\label{thm101}(\cite{Y. Caro}) Any non-complete graph with $\delta(G)\geq
n/2+\log{n}$ has $rc(G)=2$.\qed
\end{thm}

Let $G=G(n,p)$ denote, as usual, the random graph with $n$ vertices
and edge probability $p$. For a graph property $A$ and for a
function $p=p(n)$, we say that $G(n,p)$ satisfies $A$
$almost~surely$ if the probability that $G(n,p(n))$ satisfies $A$
tends to 1 as $n$ tends to infinity. We say that a function $f(n)$
is a $sharp~threshold~function$ for the property $A$ if there are
two positive constants $c$ and $C$ so that $G(n,cf(n))$ almost
surely does not satisfy $A$ and $G(n,p)$ satisfies $A$ almost surely
for all $p\geq Cf(n)$. In \cite{Y. Caro}, the authors also
determined the threshold function for a random graph to have
$rc(G(n,p))\leq 2$.

\begin{thm}\label{thm102}(\cite{Y. Caro}) $p=\sqrt{\log{n}/n}$ is a
sharp threshold function for the property $rc(G(n,p))\leq 2$.\qed
\end{thm}

In \cite{Chartrand 1}, the authors derived that the following
proposition.

\begin{pro}\label{thm103}(\cite{Chartrand 1}) $rc(G)=2$ if and only
if $src(G)=2$.\qed
\end{pro}

\begin{figure}[!hbpt]
\begin{center}
\includegraphics[scale=1.000000]{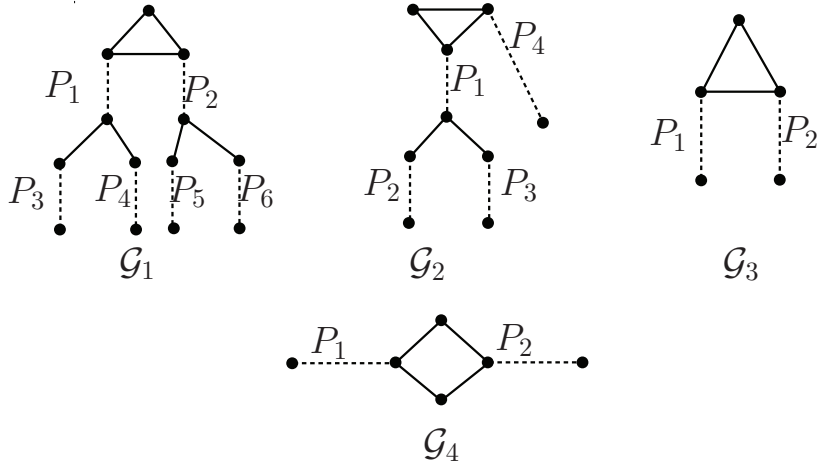}
\end{center}
\caption{The figure for the four graph classes.}\label{figure6}
\end{figure}

That is, the problem of considering graphs with $rc(G)=2$ is
equivalent to that of considering graphs with $src(G)=2$. So we aim
to investigate the graphs with large (strong) rainbow connection
numbers. In \cite{Li-Sun 1}, we investigated the graphs with large
rainbow connection numbers. In \cite{Chartrand 1}, Chartrand et al.
obtained that $rc(G)=m$ if and only if $G$ is a tree. And in
\cite{Li-Sun 1}, we proved that $rc(G)\neq m-1$, so $rc(G)\leq m-2$
if and only if $G$ is not a tree. Furthermore, we characterized the
graphs with $rc(G)=m-2$. The four graph classes shown in Figure
\ref{figure6} was useful in the following result, where the paths
$P_j$s may be trivial in each $\mathcal{G}_i$$(1\leq i\leq 4)$.
\begin{thm}\label{thm104}(\cite{Li-Sun 1}) $rc(G)\leq m-2$ if
and only if $G$ is not a tree. Furthermore, $rc(G)= m-2$ if and only
if $G$ is a 5-cycle or belongs to one of the four graph classes
shown in Figure \ref{figure6}.\qed
\end{thm}

In this paper, we continue to investigate the graphs with large
strong rainbow connection numbers with an essentially different
argument. In \cite{Chartrand 1}, Chartrand et al. obtained that
$src(G)=m$ if and only if $G$ is a tree. Furthermore, we can show
that $src(G)\neq m-1$, so $src(G)\leq m-2$ if and only if $G$ not is
a tree. We also characterize the graphs with $src(G)=m-2$ by showing
that $src(G)=m-2$ if and only if $G$ is a 5-cycle or belongs to one
of three graph classes (Theorem \ref{thm1}).

In \cite{Chartrand 1}, Chartrand et al. determined the precise
strong rainbow connection numbers for other special graph classes
including complete graphs, wheels, complete bipartite (multipartite)
graphs. However, for a general graph $G$, it is almost impossible to
give the precise value for $src(G)$, so we aim to give upper bounds
for it. In this paper, we will derive a sharp upper bound for
$src(G)$ according to the number of edge-disjoint triangles in graph
$G$, and give a necessary and sufficient condition for the equality.

We use~$V(G)$, $E(G)$ for the set of vertices and edges of $G$,
respectively. For any subset $X$ of $V(G)$, let $G[X]$ be the
subgraph induced by $X$, and $E[X]$ the edge set of $G[X]$;
similarly, for any subset $E_1$ of $E(G)$, let $G[E_1]$ be the
subgraph induced by $E_1$. Let $\mathcal {G}$ be a set of graphs,
then $V(\mathcal {G})=\bigcup_{G\in \mathcal {G}}{V(G)}$,
$E(\mathcal {G})=\bigcup_{G\in \mathcal {G}}{E(G)}$. A
$rooted~tree~T(x)$ is a tree $T$ with a specified vertex $x$, called
the $root$ of $T$. Each vertex on the path $xTv$, including the
vertex $v$ itself, is called an $ancestor$ of $v$, an ancestor of a
vertex is $proper$ if it is not the vertex itself, the immediate
proper ancestor of a vertex $v$ other than the root is its $parent$
and the vertices whose parent is $v$ are its $children$ or $son$. We
let $P_n$ and $C_n$ be the path and cycle with $n$ vertices,
respectively. $P:u_1,u_2,\cdots,u_t$ is a path, then the $u_i-u_j$
$section$ of $P$, denoted by $u_iPu_j$, is the path:
$u_i,u_{i+1},\cdots,u_j$. Similarly, for a cycle $C:
v_1,\cdots,v_t,v_1$; we define the $v_i-v_j$ section, denoted by
$v_iCv_j$ of $C$, and $C$ contains two $v_i-v_j$ sections. Note the
fact that if $P$ is a $u_1-u_t$ geodesic, then $u_iPu_j$ is also a
$u_i-u_j$ geodesic where $1\leq i,j\leq t$. We use $l(P)$ to denote
the length of path $P$. Let $[n]=\{1,\cdots,n\}$ denote the set of
the first $n$ natural numbers. For a set $S$, $|S|$ denote the
cardinality of $S$. In a graph G which has at least one cycle, the
length of a shortest cycle is called its $girth$, denoted $g(G)$. We
follow the notation and terminology of \cite{Bondy}.

\section{Basic results}

We first give a necessary condition for an edge-colored graph to be
strong rainbow connected. If $G$ contains at least two cut edges,
then for any two cut edges $e_1=u_1u_2$, $e_1=v_1v_2$, there must
exist some $1\leq i_0,j_0\leq 2$, such that any $u_{i_0}-v_{j_0}$
path must contain edge $e_1,e_2$. So we have:
\begin{obs}\label{ob1}
If $G$ is strong rainbow connected under some edge-coloring, $e_1$
and $e_2$ are any two cut edges (if exist), then $$c(e_1)\neq
c(e_2).$$\qed
\end{obs}

We need a lemma which will be useful in the argument of our result
on graphs with large strong rainbow connection numbers.

\begin{figure}[!hbpt]
\begin{center}
\includegraphics[scale=1.000000]{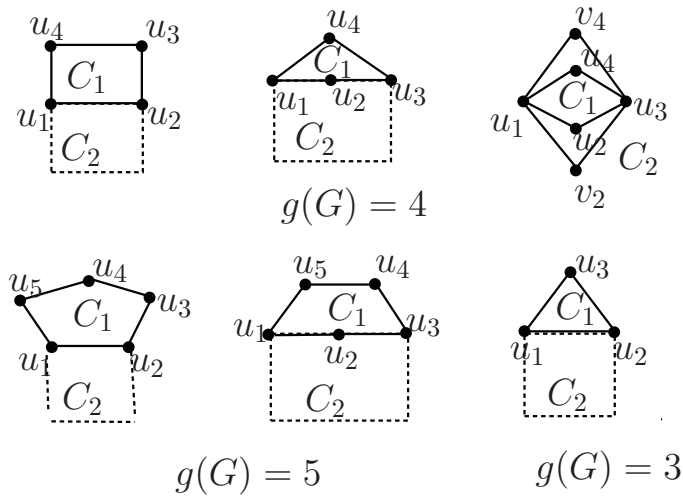}
\end{center}
\caption{Figure for Lemma \ref{lem1}.}\label{figure2}
\end{figure}

\begin{lem}\label{lem1}
$G$ is a connected graph with at least one cycle, and $3\leq g(G)\leq 5$.
Let $C_1$ be the smallest cycle of $G$, and
$C_2$ be the smallest cycle among all remaining cycles (if exist) of $G$.
If $C_1$ and $C_2$ have at least two common
vertices, then we have:

\textbf{1.} If $g(G)=3$, then $C_1$ and $C_2$ have a common edge as
shown in Figure \ref{figure2};

\textbf{2.} If $g(G)=4$, then $C_1$ and $C_2$ have a common edge, or
two common (adjacent) edges, or $C_1$ and $C_2$ are two
edge-disjoint 4-cycles, as shown in Figure \ref{figure2};

\textbf{3.} If $g(G)=5$, then $C_1$ and $C_2$ have a common edge, or
two common (adjacent) edges, as shown in Figure \ref{figure2}.

\end{lem}
\begin{pf} We only consider the case that $g(G)=5$, the remaining
two cases are similar. Let $|C_i|=k_i(i=1,2)$, and $V(C_1)=\{u_i:
1\leq i\leq 5\}$, we have $k_2\geq k_1=5$. We will consider four
cases according to the value of $|V(C_1)\cap V(C_2)|$.

\textbf{Case 1.} $|V(C_1)\cap V(C_2)|=5$, that is, $V(C_1)\subseteq
V(C_2)$. By the choice of $C_1$ and $C_2$, we have $|E(C_1)\cap
E(C_2)|\leq 4$. Without loss of generality, let $u_1u_2$ doesn't
belong to $E(C_2)$. Let $C_2'$ be the $u_1C_2u_2$ section of $C_2$
which doesn't contain $u_3,u_4,u_5$, then $C_2'$ and $u_1u_2$
produce a smaller cycle than $C_2$, a contradiction.

\textbf{Case 2.} $|V(C_1)\cap V(C_2)|=4$, say $u_i (1\leq i\leq 4)
\in V(C_2)$. In this case, we have $|E(C_1)\cap E(C_2)|\leq 3$.

\textbf{Subcase 2.1.} $|E(C_1)\cap E(C_2)|=3$, that is,
$u_1u_2,u_2u_3,u_3u_4\in E(C_2)$. Let $C_2'$ be the $u_1C_2u_4$
section of $C_2$ which doesn't contain $u_2,u_3$, then $C_2'$ and
edges $u_1u_5, u_5u_4$ produce a smaller cycle than $C_2$, a
contradiction.

\textbf{Subcase 2.2.} $|E(C_1)\cap E(C_2)|\leq 2$. Then there exists
one edge, say $u_3u_4$, which doesn't belong to $C_2$. Let $C_2'$ be
the $u_3C_2u_4$ section of $C_2$ which doesn't contain $u_1,u_2$,
then $C_2'$ and edge $u_3u_4$ produce a smaller cycle than $C_2$, a
contradiction.

\textbf{Case 3.} $|V(C_1)\cap V(C_2)|=3$.

\textbf{Subcase 3.1.} These three vertices are consecutive, say
$u_1,u_2,u_3\in  V(C_2)$. In this case, we have $|E(C_1)\cap
E(C_2)|\leq 2$. Suppose $|E(C_1)\cap E(C_2)|\leq 1$. Without loss of
generality, we assume $u_2u_3$ doesn't belong to $C_2$, and let
$C_2'$ be the $u_2C_2u_3$ section of $C_2$ that doesn't contain
$u_1$, $C_2''$ be the $u_1C_2u_3$ section of $C_2$ that doesn't
contain $u_2$. Then $C_2''$ and edges $u_1u_5,u_5u_4,u_4u_3$ will
produce a cycle $\overline{C_2}$ with $k_3=|\overline{C_2}|\geq
k_2\geq 5$, so $l(C_2'')\geq 2$, but now $C_2'$ and $u_2u_3$ produce
a smaller cycle than $C_2$, a contradiction. So we have $|E(C_1)\cap
E(C_2)|= 2$.

\textbf{Subcase 3.2.} Two of these three vertices are not
consecutive, say $u_1,u_2,u_4$. In this case, we have $|E(C_1)\cap
E(C_2)|\leq 1$. With a similar argument to \textbf{Subcase 3.1}, we
get a contradiction.

\textbf{Case 4.} $|V(C_1)\cap V(C_2)|=2$.

\textbf{Subcase 4.1.} These two vertices are adjacent in $C_1$, say
$u_1,u_2$. If $u_1u_2$ doesn't belong to $C_2$, then edge $u_1u_2$
is a chord of $C_2$ which divides $C_2$ into two parts $C_2',C_2''$.
Let $\overline{C_2'}$ ($\overline{C_2''}$) be the cycle produced by
edge $u_1u_2$ and $C_2'$ $(C_2'')$. So we have
$|\overline{C_2'}|,|\overline{C_2''}|\geq 5$, $|C_2'|,|C_2''|\geq
4$, so we have $|\overline{C_2'}|,|\overline{C_2''}|<|C_2'|+
|C_2''|=|C_2|$ and get a contradiction, so we have $u_1u_2\in C_2$.

\textbf{Subcase 4.2.} These two vertices are nonadjacent in $C_1$,
say $u_1,u_3$. Then with a similar argument to \textbf{Subcase
4.1}(instead $u_1u_2$ by $u_1u_2,u_2u_3$), we get a contradiction.

So by the above discussion, $\textbf{3}$ holds.\qed
\end{pf}

\section{Graphs with large strong rainbow connection numbers}

In this section, we will give our result on graphs with large strong
rainbow connection numbers. We first introduce three graph classes.
Let $C$ be the cycle of a unicyclic graph $G$,
$V(C)=\{v_1,\cdots,v_k\}$ and $\mathcal {T}_{G}=\{T_i:1\leq i\leq k
\}$ where $T_i$ is the unique tree containing vertex $v_i$ in
subgraph $G\backslash E(C)$. We say $T_i$ and $T_j$ are
$adjacent(nonadjacent)$ if $v_i$ and $v_j$ are adjacent(nonadjacent)
in cycle $C$. Then
let\\
$\mathcal {G}_1=\{G: G$ is a unicyclic graph, $k=3$, $\mathcal
{T}_{G}$ contains at most two nontrivial elements$\}$,\\
$\mathcal {G}_2=\{G: G$ is a unicyclic graph, $k=4$, $\mathcal
{T}_{G}$ contains two nonadjacent nontrivial
elements and each nontrivial element is a path, the remaining two elements are trivial$\}$,\\
$\mathcal {G}_3=\{G: G$ is a unicyclic graph, $k=4$, $\mathcal
{T}_{G}$ contains at most one nontrivial element and this nontrivial
element (if exists) is a path$\}$.

The following theorem is one of our main results. During its proof,
we derive that $src(G)\neq m-1$.

\begin{thm}\label{thm1}
$src(G)\leq m-2$ if and only if $G$ is not a tree, furthermore,
$src(G)=m-2$ if and only if $G$ is a 5-cycle or belongs to one of
$\mathcal {G}_i$s$(1\leq i\leq 3)$.
\end{thm}
\begin{pf} In \cite{Chartrand 1}, the authors obtained that $src(G)=m$ if
and only if $G$ is a tree, so $src(G)\leq m-1$ if and only if $G$ is
not a tree. In order to derive our conclusion, we need the following
three claims.

\textbf{Claim 1.} If $src(G)=m-1$ or $m-2$, then $3\leq g(G)\leq 5$.

\textbf{Proof of Claim 1.}
Let $C:v_1,\cdots,v_k,v_{k+1}=v_1$ be a minimum cycle of $G$ with $k=g(G)$, and
$e_i=v_iv_{i+1}$ for each $1\leq i\leq k$, we suppose that $k\geq 6$.

Now we give the cycle $C$ a strong rainbow coloring the same as
\cite{Chartrand 1}: If $k$ is even, let $k=2\ell$ for some integer
$\ell \geq 3$, $c(e_i)=i$ for $1\leq i\leq \ell$ and $c(e_i)=i-\ell$
for $\ell+1\leq i\leq k$; If $k$ is odd, let $k=2\ell+1$ for some
integer $\ell \geq 3$, $c(e_i)=i$ for $1\leq i\leq \ell+1$ and
$c(e_i)=i-\ell-1$ for $\ell+2\leq i\leq k$. We color each other edge
with a fresh color. This procedure costs $\lceil
\frac{k}{2}\rceil+(m-k)=m-(k-\lceil \frac{k}{2}\rceil)\leq m-3$
colors totally.

We will show that, with the above coloring, $G$ is strong rainbow
connected, it suffices to show that there is a rainbow $u-v$
geodesic for any two vertices $u, v$ of $G$. We first consider the
case $k=2\ell(\ell \geq 3)$. If there exists one $u-v$ geodesic $P$
which have at most one common edge with $C$, then $P$ must be a
rainbow geodesic.

So we can assume that each $u-v$ geodesic have at least two common
edges with $C$, we choose one such geodesic, say $P: u=u_1, \cdots,
v=u_t$. If there are two edges of $P$, say $e_1'$, $e_2'$, with the
same color, then they must be in $C$, too. Without loss of
generality, let $e_1'=v_1v_2$, we first consider the case that
$e_1'=v_1v_2$, and $v_1=u_{i_1}, v_2=u_{i_1+1}$ for some $1\leq
i_1\leq t$, then we must have $e_2'=v_{\ell+1}v_{\ell+2}$ where
$v_{\ell+1}=u_{j_1}$, $v_{\ell+2}=u_{j_1+1}$ for some $i_1+1\leq
j_1\leq t$ or $v_{\ell+2}=u_{j_2}$, $v_{\ell+1}=u_{j_2+1}$ for some
$i_1+1\leq j_2\leq t$. If $v_{\ell+1}=u_{j_1}$,
$v_{\ell+2}=u_{j_1+1}$ for some $i_1+1\leq j_1\leq t$ (For example,
see graph $(\alpha)$ of Figure \ref{figure4} where $\ell=4$, the
color of each edge is shown), then the section $v_2Pv_{\ell+1}$ of
$P$ is a $v_2-v_{\ell+1}$ geodesic, so it is not longer than the
section $C':v_2,v_3,\cdots,v_{\ell+1}$ of $C$, then the length of
$v_2Pv_{\ell+1}$, $l(v_2Pv_{\ell+1})\leq \ell-1$, is smaller than
the length of the section $C'':v_2, v_1, v_{k}, \cdots, v_{\ell+1}$
of $C$. So the sections $v_2Pv_{\ell+1}$ and $C'$ will produce a
smaller cycle than $C$(this produces a contradiction), or
$v_2Pv_{\ell+1}$ is the same as $C'$(but in this case, the section
$C''':v_1, v_k,\cdots,v_{\ell+2}$ of $C$ is shorter than
$v_1Pv_{\ell+2}$ which now is a $v_1-v_{\ell+2}$ geodesic, this also
produces a contradiction). If $v_{\ell+2}=u_{j_2}$,
$v_{\ell+1}=u_{j_2+1}$ for some $i_1+1\leq j_2\leq t$ (For example,
see graph $(\beta)$ of Figure \ref{figure4} where $\ell=4$, the
color of each edge is shown), then the section $v_1Pv_{\ell+2}$ of
$P$ is a $v_1-v_{\ell+2}$ geodesic, so it is not longer than the
length of the section
$\overline{C'}:v_1,v_k,v_{k-1},\cdots,v_{\ell+2}$ of $C$ and its
length, $l(v_1Pv_{\ell+2})\leq \ell-1$, is smaller than that of the
section $\overline{C''}:v_1,v_2,\cdots,v_{\ell+2}$ of $C$. So the
sections $v_1Pv_{\ell+2}$ and $\overline{C'}$ will produce a smaller
cycle than $C$, this also produces a contradiction.

The remaining two subcases correspond to the case that
$v_1=u_{i_1+1}$, $v_2=u_{i_1}$(see graphs $(\gamma)$ and $(\omega)$
in Figure \ref{figure4} for the case of $\ell=4$), and with a
similar argument, a contradiction will be produced. So $P$ is
rainbow.
\begin{figure}[!hbpt]
\begin{center}
\includegraphics[scale=1.000000]{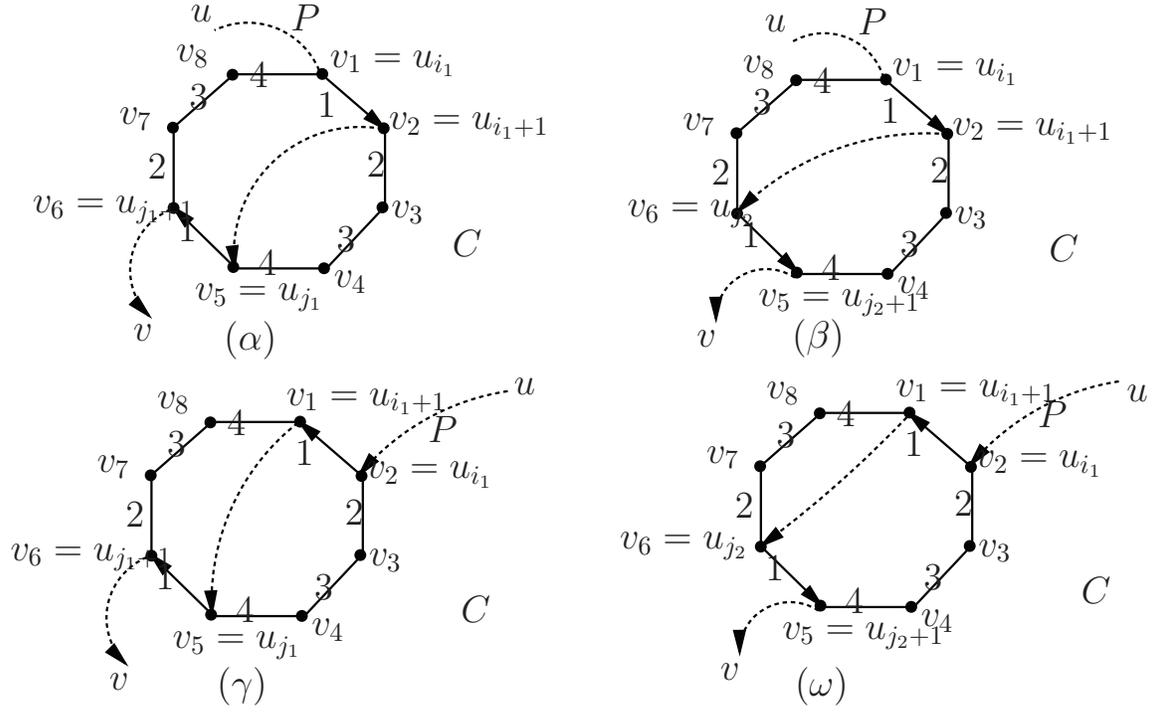}
\end{center}
\caption{Graphs for the example with $\ell=4$.}\label{figure4}
\end{figure}

The case that $k=2\ell+1(\ell\geq 3)$ is similar. So $G$ does not
contain a cycle of length larger than 5.

Note that during the proof of \textbf{Claim 1}, we use the following
\textbf{technique}: we first choose a smallest cycle $C$ of a graph
$G$, then give it a strong rainbow coloring the same as
\cite{Chartrand 1}, and give a fresh color to any other edge. Then
for any $u-v$ geodesic $P$, we derive that either one section of $P$
is the same as one section of $C$ and then find a shorter path than
the geodesic, or one section of $P$ and one section of $C$ produce a
smaller cycle than $C$, each of these two cases will produce a
contradiction.  This technique will be useful in the sequel.

Next we will show that $G$ is a unicyclic graph under the condition
that $src(G)=m-1$ or $m-2$.

\textbf{Claim 2.} If $src(G)=m-1$ or $m-2$, then $G$ is a unicyclic graph.

\textbf{Proof of Claim 2.} Suppose $G$ contains at least two cycles,
let $C_1$ be the smallest cycle of $G$ and $C_2$ be the smallest one
among all the remaining cycles in $G$, that is, $C_2$ is the
smallest cycle with the exception of cycle $C_1$. Let
$|C_i|=k_i(i=1,2)$, so by the above discussion, we have $3\leq
k_1\leq 5$ and $k_2\geq k_1$. We will consider two cases according
to the value of $|V(C_1)\cap V(C_2)|$.

\textbf{Case 1.}$|V(C_1)\cap V(C_2)|\leq 1$, that is, $C_1$ and
$C_2$ have at most one common vertex. There are three subcases:

\textbf{Subcase 1.1.} $k_1=3$, that is, $C_1$ is a triangle.
The following fact is trivial and will be useful:

\textbf{Fact 1.} For any two vertices $u$, $v$ and a triangle $T$,
any $u-v$ geodesic $P$ contains at most one edge of
$T$.

We first give cycle $C_2$ a strong rainbow coloring using $\lceil
\frac{k_2}{2} \rceil$ colors the same as \cite{Chartrand 1}; then
give a fresh color to $C_1$, that is, edges of $C_1$ receive the
same color; for the remaining edges, we give each of them a fresh
color. With a similar procedure (technique) to that of \textbf{Claim
1} and by \textbf{Fact 1}, we can show that the above coloring is
strong rainbow, as this costs $1+\lceil
\frac{k_2}{2}\rceil+(m-k_2-3)$ colors totally, we have $src(G)\leq
1+\lceil \frac{k_2}{2}\rceil+(m-k_2-3)= (m-2)-(k_2- \lceil
\frac{k_2}{2} \rceil) \leq m-3$, a contradiction.

\textbf{Subcase 1.2.} $k_1=4$, that is, $C_1$ is a 4-cycle. With a
similar (and a little simpler) argument to that of \textbf{Claim 1},
we can give the following fact,

\textbf{Fact 2.} For any two vertices $u$, $v$, any $u-v$ geodesic
$P$ contains at most one edge or two (adjacent) edges of $C_1$.

We now give a brief proof to \textbf{Fact 2}: Suppose it doesn't
holds, that is, there exist a geodesic $P: a_1,\cdots,a_t$ for two
vertices $u,v$ which contains two nonadjacent edges, say $u_1u_2,
u_3u_4$. Without loss of generality, we let $u_1=a_{i_1},
u_2=a_{i_2}, u_3=a_{i_3}, u_4=a_{i_4}$ where $\max\{i_1,i_2 \} <\min
\{i_3,i_4\}$. We only consider the case that $i_1<i_2, i_3<i_4$, the
remaining three cases are similar. Then the section
$a_{i_1}Pa_{i_4}$ of $P$ is a $u_1-u_4$ geodesic whose length is at
least three, but the edge $u_1u_4$ is a $u_1-u_4$ path which is
shorter than it, this produces a contradiction, so the fact holds.

We first give cycle $C_2$ a strong rainbow coloring using $\lceil
\frac{k_2}{2} \rceil$ colors the same as \cite{Chartrand 1}; then
give two fresh colors to $C_1$ in the same way; for the remaining
edges, we give each of them a fresh color. With a similar procedure
(technique) to that of \textbf{Claim 1} and by \textbf{Fact 2}, we
can show that the above coloring is strong rainbow, as this costs
$2+\lceil \frac{k_2}{2}\rceil+(m-k_2-4)$ colors totally, we have
$src(G)\leq 2+\lceil \frac{k_2}{2}\rceil+(m-k_2-4)=
(m-2)-(k_2-\lceil \frac{k_2}{2} \rceil)\leq m-3$, a contradiction.

\textbf{Subcase 1.3.} $k_1=5$, that is, $C_1$ is a 5-cycle. With a
similar argument to that of \textbf{Fact 2}, we can give the
following fact,

\textbf{Fact 3.} For any two vertices $u$, $v$, any $u-v$ geodesic
$P$ contains at most one edge or two (adjacent) edges of $C_1$.

We first give cycle $C_2$ a strong rainbow coloring using $\lceil
\frac{k_2}{2} \rceil$ colors the same as \cite{Chartrand 1}; then
give three fresh colors to $C_1$ in the same way; for the remaining
edges, we give each of them a fresh color. With a similar procedure
(technique) to that of \textbf{Claim 1} and by \textbf{Fact 3}, we
can show that the above coloring is strong rainbow, as this costs
$3+\lceil \frac{k_2}{2}\rceil+(m-k_2-5)$ colors totally, we have
$src(G)\leq 3+\lceil \frac{k_2}{2}\rceil+(m-k_2-5)=
(m-2)-(k_2-\lceil \frac{k_2}{2} \rceil)\leq m-3$, a contradiction.

Note that for each above subcase, by each corresponding fact, the
cycle produced during the procedure while we use the technique of
\textbf{Claim 1} cannot be the cycle $C_1$ and must be smaller than
$C_2$, then a contradiction will be produced.

So if $src(G)=m-1$ or $m-2$, \textbf{Case 1} doesn't hold, we now
consider the next case:

\textbf{Case 2.} $|V(C_1)\cap V(C_2)|\geq 2$, that is, $C_1$ and
$C_2$ have at least two common vertices. Similarly, we need to
consider three subcases:
\begin{figure}[!hbpt]
\begin{center}
\includegraphics[scale=1.000000]{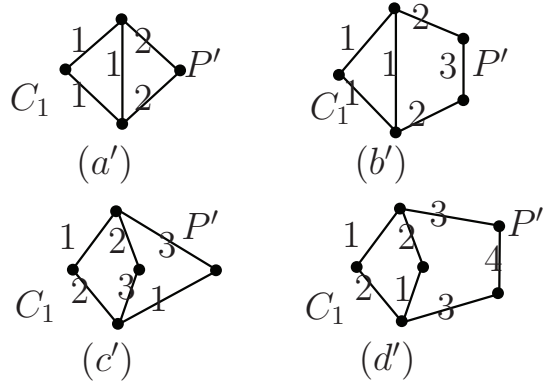}
\end{center}
\caption{Graphs for Case 2 of Claim 2.}\label{figure3}
\end{figure}

\textbf{Subcase 2.1.} $k_1=3$, that is, $C_1$ is a triangle. By
Lemma \ref{lem1}, $C_1$ and $C_2$ have a common edge as shown in
Figure \ref{figure2}. Let $V(C_2)=\{v_i:1\leq i\leq k_2\}$ and
$v_{k_2+1}=v_1$, where $v_1=u_1,v_2=u_2$. Let $P'$ be the subpath of
$C_2$ that doesn't contain edge $v_1v_2$. We now give $G$ an
edge-coloring as follows:

For the case $l(P')=2,3$, we first color edges of $C_1\cup C_2$ as
shown in Figure \ref{figure3} (graphs $a'$ and $b'$); then give each
other edge of $G$ a fresh color. This procedure costs $m-3$ colors
totally. Then with a similar argument to \textbf{Fact 2}, we can
show that any geodesic cannot contain two edges with the same color,
so $src(G)\leq m-3$.

For the remaining case, that is, $l(P')\geq 4$ and $k_2\geq 5$. We
first give cycle $C_1$ a color, say $a$, that is, three edges of
$C_1$ receive the same color. Then in $C_2$, if $k_2=2\ell$ for some
$\ell \geq 2$, then let $c(v_2v_3)=c(v_{\ell+2}v_{\ell+3})$ be a new
color, say $b$; if $k_2=2\ell+1$ for some $\ell \geq 2$, then let
$c(v_2v_3)=c(v_{\ell+3}v_{\ell+4})$ be a new color, say $b$. For the
remaining edges, we give each of them a fresh color. This procedure
costs $m-3$ colors totally. For any two vertices $u,v$, $P$ is a
$u-v$ geodesic, by \textbf{Fact 1}, $P$ cannot contain two edges
with color $a$; for the two edges with color $b$, with a similar
argument to that of \textbf{Claim 1} (Note that now, by \textbf{Fact
1}, the cycle produced during the procedure cannot be $C_1$ and must
be shorter than $C_2$, then a contradiction will be produced), we
can show $P$ contains at most one of them. So $P$ is rainbow and
$src(G)\leq m-3$.

\textbf{Subcase 2.2.} $k_1=4$, that is, $C_1$ is a 4-cycle. By Lemma
\ref{lem1}, $C_1$ and $C_2$ have a common edge, or two common
(adjacent) edges, or $C_1$ and $C_2$ are two edge-disjoint 4-cycles,
as shown in Figure \ref{figure2}.

\textbf{Subsubcase 2.2.1.} If $C_1$ and $C_2$ have a common edge,
say $u_1u_2$ (see the left one of the three graphs with $g(G)=4$ in
Figure \ref{figure2}). We let $V(C_2)=\{v_i:1\leq i\leq k_2\}$,
where $v_1=u_1,v_2=u_2$. We let $c(v_2v_3)=c(u_4v_1)=a$,
$c(v_2u_3)=c(v_1v_{k_2})=b$, $c(v_1v_2)=c(u_3u_4)=c$. For the
remaining edges, we give each of them a fresh color. This procedure
costs $m-3$ colors totally. For any two vertices $u,v$, $P$ is a
$u-v$ geodesic, then by \textbf{Fact 2}, $P$ contains at most one of
the two edges with color $c$; for the two edges with color $a(b)$,
with a similar argument to that of \textbf{Fact 2}, we can show that
there exists one $u-v$ geodesic which contains at most one of them.
So we have $src(G)\leq m-3$.

\textbf{Subsubcase 2.2.2.} If $C_1$ and $C_2$ have two (adjacent)
common edges, say $u_1u_2,u_2u_3$ (see the middle one of the three
graphs with $g(G)=4$ in Figure \ref{figure2}). We let
$V(C_2)=\{v_i:1\leq i\leq k_2\}$, where $v_1=u_1,v_2=u_2,v_3=u_3$.
Let $P'$ be the subpath  of $C_2$ which doesn't contain edges
$u_1u_2,u_2u_3$.

For the case $l(P')=2,3$, we first color edges of $C_1\cup C_2$ as
shown in Figure \ref{figure3} (graphs $c'$ and $d'$); then give each
other edge of $G$ a fresh color. This procedure costs $m-3$ colors
totally. Then with a similar argument to \textbf{Fact 2}, we can
show that any geodesic cannot contain two edges with the same color,
so we have $src(G)\leq m-3$.

For the case $l(P')\geq 4$, that is $k_2\geq 6$. We let
$c(u_4v_1)=c(v_3v_4)=a$, $c(v_1v_2)=c(v_3u_4)=b$; for edge $v_2v_3$,
we give a similar treatment to that of \textbf{Subcase 2.1} and let
$c(v_2v_3)=c$; we then give each other edge of $G$ a fresh color.
This procedure costs $m-3$ colors totally. For any two vertices
$u,v$, $P$ is a $u-v$ geodesic, then by \textbf{Fact 2}, $P$
contains at most one of the two edges with color $b$; for the two
edges with color $a$, with a similar argument to that of
\textbf{Fact 2}, we can show that there exists one $u-v$ geodesic
which contains at most one of them. With a similar argument to that
of \textbf{Claim 1}(Note that now, by \textbf{Fact 2}, the cycle
produced during the procedure cannot be $C_1$ and must be shorter
than $C_2$, then a contradiction will be produced), we can show any
geodesic contains at most one edge with color $c$. So we have
$src(G)\leq m-3$.

\textbf{Subsubcase 2.2.3.} The remaining case, the right graph of
the three graphs with $g(G)=4$ in Figure \ref{figure2}. We let
$c(u_1u_2)=c(u_3u_4)=a, c(u_2u_3)=c(u_1u_4)=b,
c(u_1v_2)=c(u_3v_4)=c, c(v_2u_3)=c(u_1v_4)=d,$ where $a,b,c,d$ are
four distinct colors; for the remaining edges, we give each of them
a fresh color. This procedure costs $m-4$ colors totally. As now
both $C_1$ and $C_2$ are the smallest cycle of $G$, by \textbf{Fact
2}, any geodesic contains at most one of the two edges with the same
color, so $src(G)\leq m-4$.

\textbf{Subcase 2.3.} $k_1=5$, that is, $C_1$ is a 5-cycle. By Lemma
\ref{lem1}, $C_1$ and $C_2$ have a common edge, or two common
(adjacent) edges, as shown in Figure \ref{figure2}. The following
discussion will use \textbf{Fact 3}.

\textbf{Subsubcase 2.3.1.} If $C_1$ and $C_2$ have a common edge,
say $u_1u_2$ (see the left one of the two graphs with $g(G)=5$ in
Figure \ref{figure2}). We let $V(C_2)=\{v_i:1\leq i\leq k_2\}$,
where $v_1=u_1,v_2=u_2$. We let $c(u_4u_5)=c(v_2v_3)=a$,
$c(v_1u_5)=c(v_2u_3)=b$, and $c(v_1v_2)=c(u_3u_4)=c$; for the
remaining edges, we give each of them a fresh color. This procedure
costs $m-3$ colors totally. With a similar argument to above, we can
show that $src(G)\leq m-3$.

\textbf{Subsubcase 2.3.2.} If $C_1$ and $C_2$ have two common
(adjacent) edges, say $u_1u_2,u_2u_3$ (see the right one of the two
graphs with $g(G)=5$ in Figure \ref{figure2}). We let
$c(v_1u_5)=c(v_3v_4)=a$, $c(v_1v_2)=c(v_3u_4)=b$, and
$c(v_2v_3)=c(u_4u_5)=c$; for the remaining edges, we give each of
them a fresh color. This procedure costs $m-3$ colors totally. With
a similar argument to above, we can show that $src(G)\leq m-3$.

With the above discussion, \textbf{Claim 2} holds. To complete our
proof, we still need the following claim.

\textbf{Claim 3.} $src(G)\neq m-1$. Furthermore, $src(G)=m-2$ if and
only if $G$ is a 5-cycle or belongs to one of $\mathcal
{G}_i$s$(1\leq i\leq 3)$.

\textbf{Proof of Claim 3.} Let $G$ be a unicyclic graph and $C$ be
its cycle, $|C|=k$. We consider three cases according to the length
$k$ of cycle $C$.

\textbf{Case 1.} $k=3$.

\textbf{Subcase 1.1.} All $T_i$s are nontrivial. We first give each
edge of $G\backslash E(C)$ a fresh color, then let $c(v_1v_2)\in
c(T_3)$, $c(v_2v_3)\in c(T_1)$, $c(v_1v_3)\in c(T_2)$, it is easy to
show, with this coloring, $G$ is strong rainbow connected, so
$src(G)\leq m-3$ in this case.

\textbf{Subcase 1.2.} At most two $T_i$s are nontrivial, that is,
$G\in \mathcal{G}_1$. We first consider the case that there are
exactly two $T_i$s which are nontrivial, say $T_1$ and $T_2$. We
first give each edge of $G\backslash E(C)$ a fresh color, then let
$c(v_1v_2)=c(v_2v_3)=c(v_1v_3)$, it is easy to show, with this
coloring, $G$ is strong rainbow connected, so now $src(G)\leq m-2$.
On the other hand, by Observation \ref{ob1} and the definition of
rainbow geodesic, we know that $c(T_1)\cap c(T_2)=\emptyset$ and
$c(v_1v_2)$ doesn't belong to $c(T_1)\cup c(T_2)$. So we have
$src(G)= m-2$ in this case. With a similar argument, we can derive
$src(G)= m-2$ for the case that at most one $T_i$ is nontrivial. So
$src(G)= m-2$ if $G\in \mathcal{G}_1$.

\textbf{Case 2.} $k=4$.

\textbf{Subcase 2.1.} There are at least three nontrivial $T_i$s,
say $T_1,T_3,T_4$. We first give each edge of $G\backslash E(C)$ a
fresh color, then let $c(v_1v_2)\in c(T_3)$, $c(v_3v_4)\in c(T_1)$,
$c(v_1v_4)\in c(T_2)$ and we give edge $v_2v_3$ a fresh color. It is
easy to show, with this coloring, $G$ is strong rainbow connected,
so $src(G)\leq m-3$ in this case.

\textbf{Subcase 2.2.} There are exactly two nontrivial $T_i$s, say
$T_{i_1}$ and $T_{i_2}$.

\textbf{Subsubcase 2.2.1.} $T_{i_1}$ and $T_{i_2}$ are adjacent, say
$T_1$ and $T_2$. We first give each edge of $G\backslash E(C)$ a
fresh color, then let $c(v_2v_3)\in c(T_1)$, $c(v_1v_4)\in c(T_2)$
and we color edges $v_1v_2$ and $v_3v_4$ with the same new color. It
is easy to show, with this coloring, $G$ is strong rainbow
connected, so $src(G)\leq m-3$ in this case.

\begin{figure}[!hbpt]
\begin{center}
\includegraphics[scale=1.000000]{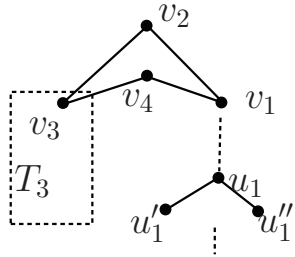}
\end{center}
\caption{Graph for Subsubcase 2.2.2 of Claim 3.}\label{figure5}
\end{figure}
\textbf{Subsubcase 2.2.2.} $T_{i_1}$ and $T_{i_2}$ are nonadjacent,
say $T_1$ and $T_3$. We can consider $T_i$ as rooted tree with root
$v_i$ $(i=1,3)$. If there exists some $T_i$, say $T_1$, that
contains a vertex, say $u_1$, with at least two sons, say
$u_1',u_1''$(see Figure \ref{figure5}). We first color each edge of
$\bigcup_{i=1,3} T_i\cup \{v_1v_2\}$ with a distinct color, this
costs $m-3$ colors, then we let $c(v_1v_4)=c(v_1v_2),
c(v_2v_3)=c(u_1u_1'), c(v_3v_4)=c(u_1u_1'')$. It is easy to show
that this coloring is strong rainbow and we have $src(G)\leq m-3$ in
this case. If $G\in \mathcal{G}_2$, we first give each edge of
$G\backslash E(C)$ a fresh color, then let $c(v_1v_2)=c(v_3v_4)=a$
and $c(v_2v_3)=c(v_1v_4)=b$ where $a$ and $b$ are two new colors. It
is easy to show, with this coloring, $G$ is strong rainbow
connected, so $src(G)\leq m-2$ in this case. On the other hand,
$src(G)\geq m-2=diam(G)$. So $src(G)= m-2$ if $G\in \mathcal{G}_2$.

\textbf{Subcase 2.3.} There are at most one nontrivial $T_i$. Then
with a similar argument to \textbf{Subsubcase 2.2.2}, we can derive
that $src(G)= m-2$ if $G\in \mathcal{G}_3$.

\textbf{Case 3.} $k=5$.

If there are at least one nontrivial $T_i$, say $T_1$, then we give
each edge of $G\backslash E(C)$ a fresh color, let $v_3v_4\in
c(T_1)$, $c(v_1v_2)=c(v_4v_5)=a$ and $c(v_2v_3)=c(v_1v_5)=b$ where
$a$ and $b$ are two new colors. It is easy to show, with this
coloring, $G$ is strong rainbow connected, so now we have
$src(G)\leq m-3$. On the other hand, we know $src(G)=m-2=3$ if
$G\cong C_5$ from \cite{Chartrand 1}.

By \textbf{Claim 1} and \textbf{Claim 2}, we have that if
$src(G)=m-1$ or $m-2$, then $G$ is a unicyclic graph with the cycle
of length at most 5. By the discussion from \textbf{Case 1} to
\textbf{Case 3} of \textbf{Claim 3}, we know that if $G$ is a
unicyclic graph with the cycle of length at most 5, then $src(G)\neq
m-1$. So $src(G)\neq m-1$ for any graph $G$ and \textbf{Claim 3}
holds. By our three claims, our theorem holds. \qed
\end{pf}

\section{Upper bound for $src(G)$ according to edge-disjoint triangles}
In this section, we give an upper bound for $src(G)$ according to
their edge-disjoint triangles in graph $G$.

Recall that a $block$ of a connected graph $G$ is a maximal
connected subgraph without a cut vertex. Thus, every block of graph
$G$ is either a maximal 2-connected subgraph or a bridge (cut edge).
We now introduce a new graph class. For a connected graph $G$, we
say $G\in \overline{\mathcal{G}}_t$, if it satisfies the following
conditions:

\textbf{$C_1$.} Each block of $G$ is a bridge or a triangle;

\textbf{$C_2$.} $G$ contains exactly $t$ triangles;

\textbf{$C_3$.} Each triangle contains at least one vertex of degree two in $G$.

\begin{figure}[!hbpt]
\begin{center}
\includegraphics[scale=1.000000]{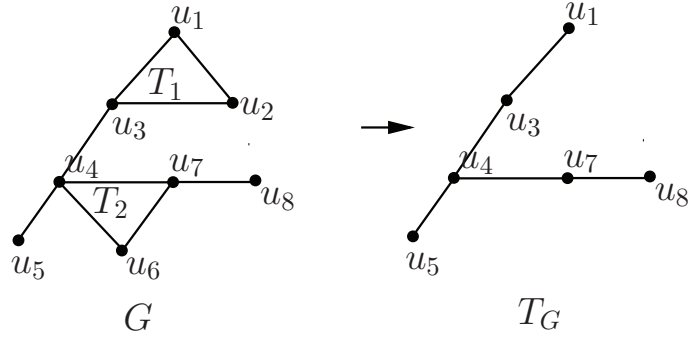}
\end{center}
\caption{An example of $G\in \overline{\mathcal{G}}_t$ with
$t=2$.}\label{figure1}
\end{figure}

By the definition, each graph $G\in \overline{\mathcal{G}}_t$ is
formed by (edge-disjoint) triangles and paths (may be trivial),
these triangles and paths fit together in a treelike structure, and
$G$ contains no cycles but the $t$ (edge-disjoint) triangles. For
example, see Figure \ref{figure1}, here $t=2$, $u_1$, $u_2$, $u_6$
are vertices of degree 2 in $G$. If a tree is obtained from a graph
$G\in \overline{\mathcal{G}}_t$ by deleting one vertex of degree 2
for each triangle, then we call this tree is a $D_2$-$tree$ of $G$,
denoted $T_{G}$. For example, in Figure \ref{figure1}, $T_{G}$ is a
$D_2$-tree of $G$. Clearly, the $D_2$-tree is not unique, since in
this example, we can obtain another $D_2$-tree by deleting vertex
$u_1$ instead of $u_2$. On the other hand, we can say any element of
$\overline{\mathcal{G}}_t$ can be obtained from a tree by adding $t$
new vertices of degree 2. It is easy to show that number of edges of
$T_{G}$ is $m-2t$ where $m$ is the number of edges of $G$.

\begin{thm}\label{thm2}
$G$ is a graph with $m$ edges and $t$ edge-disjoint triangles, then
$$src(G)\leq m-2t,$$ the equality holds if and only if $G\in
\overline{\mathcal{G}}_t$.

\end{thm}
\begin{pf}
Let $\mathcal{T}=\{ T_i: 1\leq i\leq t\}$ be a set of $t$
edge-disjoint triangles in $G$.

We color each triangle with a fresh color (this means that three
edges of each triangle receive the same color), then we give each
other edge a fresh color. For any two vertices $u,v$ of $G$, let $P$
be any $u-v$ geodesic, then $P$ contains at most one edge from each
triangle by \textbf{Fact 1}, so $P$ is rainbow under the above
coloring. As this procedure costs $m-2t$ colors totally, we have
$src(G)\leq m-2t$.

\textbf{Claim 1.} If the equality holds, then for any set
$\mathcal{T}$ of edge-disjoint triangles of $G$, we have
$|\mathcal{T}|\leq t$.

\textbf{Proof of Claim 1.} We suppose there is a set $\mathcal{T}'$
of $t'$ edge-disjoint triangles in $G$ with $t'> t$, then with a
similar procedure, we have $src(G)\leq m-2t'<m-2t$, a contradiction.

\textbf{Claim 2.} If the equality holds, then $G$ contains no cycle
but the above $t$ (edge-disjoint) triangles.

\textbf{Proof of Claim 2.} We suppose that there are at least one
cycles distinct with the above $t$ triangles. Let $\mathcal{C}$ be
the set of these cycles and $C_1$ be the smallest element of
$\mathcal{C}$ with $|C_1|=k$. We will consider two cases:

\textbf{Case 1.} $E(C_1)\cap E(\mathcal{T})=\emptyset$, that is,
$C_1$ is edge-disjoint with each of the above $t$ triangles. With a
similar argument to Lemma \ref{lem1}, we know $C_1$ has at most one
common vertex with each of the above $t$ triangles. In this case
$k\geq 4$ by \textbf{Claim 1}. We give $G$ an edge coloring as
follows: we first color edges of cycle $C_1$ the same as
\cite{Chartrand 1} (this is shown in the proof of \textbf{Claim 1}
of Theorem \ref{thm1}); then we color each triangle with a fresh
color; for the remaining edges, we give each one a fresh color.
Recall the fact that any geodesic contains at most one edge from
each triangle and with a similar procedure to the proof of
\textbf{Claim 1} of Theorem \ref{thm1}, we know the above coloring
is strong rainbow, as this procedure costs $\lceil
\frac{k}{2}\rceil+t+(m-k-3t)=(m-2t)+(\lceil
\frac{k}{2}\rceil-k)<m-2t$, we have $src(G)<m-2t$, this produces a
contradiction.

\textbf{Case 2.} $E(C_1)\cap E(\mathcal{T})\neq \emptyset$, that is,
$C_1$ have common edges with the above $t$ triangles, in this case
$k\geq 3$. By the choice of $C_1$, we know $|E(C_1)\cap E(T_i)|\leq
1$ for each $1\leq i\leq t$. We will consider two subcases according
to the parity of $k$.

\textbf{Subcase 2.1.} $k=2\ell$ for some $\ell \geq 2$. For example,
see graph $(\alpha)$ of Figure \ref{figure7}, here
$\mathcal{T}=\{T_1,T_2,T_3\}$, $V(C_1)=\{u_i:1\leq i\leq 6\}$,
$E(C_1)\cap E(T_1)=\{u_1u_2\}$, $E(C_1)\cap E(T_2)=\{u_4u_5\}$.
Without loss of generality, we assume that there exists a triangle,
say $T_1$, which contains edge $u_1u_2$ and let
$V(T_1)=\{u_1,u_2,w_1\}$, $G'=G\backslash E(T_1)$. If there exists
some triangle, say $T_2$, which contains edge
$u_{\ell+1}u_{\ell+2}$, we let
$V(T_2)=\{u_{\ell+1},u_{\ell+2},w_2\}$.

\begin{figure}[!hbpt]
\begin{center}
\includegraphics[scale=1.000000]{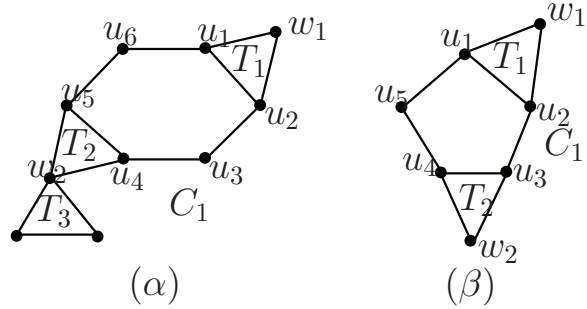}
\end{center}
\caption{Graphs of two examples in Theorem
\ref{thm2}.}\label{figure7}
\end{figure}

We first consider the case that $\ell=2$, see Figure \ref{figure8},
we first give each triangle of $G'$ a fresh color; for the remaining
edges of $G'$, we give each of them a fresh color; for edges of
$T_1$, let $c(u_1w_1)=c(u_2u_3)$, $c(u_2w_1)=c(u_1u_4)$,
$c(u_1u_2)=c(u_3u_4)$. Then with a similar argument to that of
\textbf{Fact 2}, we can show that there is a $u-v$ geodesic which
contains at most one edge from any two edges with the same color for
$u,v\in G$, so the above coloring is strong rainbow. As this
procedure costs $m-2t-1 <m-2t$ colors totally, we have
$src(G)<m-2t$, a contradiction.

We next consider the case that $\ell \geq 3$. Let $G''=G\backslash
(E(T_1)\cup E(T_2))$. We give $G$ an edge-coloring as follows: We
first give each triangle of $G''$ a fresh color; then give a fresh
color to each of the remaining edges of $G''$; for the edges of
$T_1$ and $T_2$, let $c(u_1w_1)=c(u_2u_3)=a$,
$c(u_2w_1)=c(u_1u_k)=b$, $c(u_1u_2)=c(u_{\ell+1}u_{\ell+2})=c$,
$c(w_2u_{\ell+1})=c(u_{\ell+2}u_{\ell+3})=d$,
$c(w_2u_{\ell+2})=c(u_{\ell}u_{\ell+1})=e$ where $a,b,c,d,e$ are
five new colors. Then with a similar argument to that of
\textbf{Fact 2}, we can show that there is a $u-v$ geodesic which
contains at most one edge from any two edges with the same color for
$u,v\in G$, so the above coloring is strong rainbow. As this
procedure costs $m-2t-1 <m-2t$ colors totally, we have
$src(G)<m-2t$, a contradiction.

\textbf{Subcase 2.2.} $k=2\ell+1$ for some $\ell \geq 1$.

We first consider the case that $\ell \geq 2$. For example, see
graph $(\beta)$ of Figure \ref{figure7}, here
$\mathcal{T}=\{T_1,T_2\}$, $V(C_1)=\{u_i:1\leq i\leq 5\}$,
$E(C_1)\cap E(T_1)=\{u_1u_2\}$, $E(C_1)\cap E(T_2)=\{u_3u_4\}$.
Without loss of generality, we assume that there exists a triangle,
say $T_1$, which contains edge $u_1u_2$ and let
$V(T_1)=\{u_1,u_2,w_1\}$. If there exists some triangle, say $T_2$,
which contains edge $u_{\ell+1}u_{\ell+2}$, we let
$V(T_2)=\{u_{\ell+1},u_{\ell+2},w_2\}$ and $G'=G\backslash
(E(T_1)\cup E(T_2))$.

We give $G$ an edge-coloring as follows: We first give each triangle
of $G'$ a fresh color; then give a fresh color to each of the
remaining edges of $G'$; for the edges of $T_1$ and $T_2$, let
$c(u_1w_1)=c(u_2u_3)$, $c(u_2w_1)=c(u_1u_k)$,
$c(u_{\ell+1}w_2)=c(u_{\ell+2}u_{\ell+3})$ and let
$c(u_1u_2)=c(u_{\ell+1}u_{\ell+2})=c(w_2u_{\ell+2})$ be a fresh
color. With a similar procedure to the proof of \textbf{Fact 1} and
\textbf{Claim 1} of Theorem \ref{thm1}, we can show that $G$ is
strong rainbow connected, and so $src(G)\leq
(t-1)+(m-3t)=(m-2t)-1<m-2t$, this produces a contradiction.

For the case that $\ell =1$, that is, $C_1$ is a triangle. See
Figure \ref{figure8}, we color the three edges (if exist) with color
1, these edges are shown in the figure; the remaining edges of these
three triangles (if exist) all receive color 2; each of other
triangles receive a fresh color; for the remaining edges, we give
each one a fresh color. It is easy to show that the above coloring
is strong rainbow, so we have $src(G)< m-2t$ in this case, a
contradiction. So the claim holds.

\begin{figure}[!hbpt]
\begin{center}
\includegraphics[scale=1.000000]{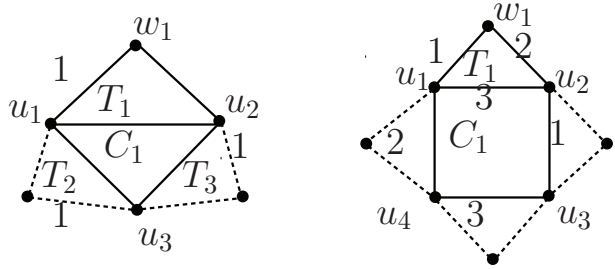}
\end{center}
\caption{Edge coloring for the case that $C_1$ is a triangle and the
case that $C_1$ a 4-cycle in Theorem \ref{thm2}.}\label{figure8}
\end{figure}

\textbf{Claim 3.} If the equality holds, then $G\in \overline{\mathcal{G}}_t$.

\textbf{Proof of Claim 3.} If the equality holds, to prove that
$G\in \overline{\mathcal{G}}_t$, it suffices to show that each
triangle contains at most one vertex of degree 2 in $G$. Suppose it
doesn't holds, without loss of generality, let $T_1$ be the triangle
with $deg_G(v_i)\geq 3$, where $V(T_1)=\{v_i: 1\leq i\leq 3 \}$.  By
\textbf{Claim 2}, it is easy to show that $E(T_1)$ is an edge-cut of
$G$, let $H_i$ be the subgraph of $G\backslash E(T_1)$ containing
vertex $v_i$ $(1\leq i\leq 3)$, by the assumption of $T_1$, we know
each $H_i$ is nontrivial. We now give $G$ an edge-coloring: for the
$t-1$ (edge-disjoint) triangles of $G\backslash E(T_1)$, we give
each of them a fresh color; for the remaining edges of $G\backslash
E(T_1)$ (by \textbf{Claim 2}, each of them must be a cut edge), we
give each of them a fresh color; for the edges of $E(T_1)$, let
$c(v_1v_3)\in c(H_2)$, $c(v_1v_2)\in c(H_3)$, $c(v_2v_3)\in c(H_1)$.
It is easy to show, with the above coloring, $G$ is strong rainbow
connected, and we have $src(G)<m-2t$, a contradiction, so the claim
holds.

\textbf{Claim 4.} If $G\in \overline{\mathcal{G}}_t$, then the equality holds.

\textbf{Proof of Claim 4.} Let $T_G$ be a $D_2$-tree of $G$, the
result clearly holds for the case $|E(T_G)|=1$. So now we assume
$|E(T_G)|\geq 2$. We will show, for any strong rainbow coloring of
$G$, $c(e_1)\neq c_(e_2)$ where $e_1,e_2 \in T_G$, that is, each
edge of $T_G$ receive a distinct color, so edges of $T_G$ cost
$m-2t$ colors totally, recall that $|E(T_G)|=m-2t$, then $src(G)\geq
m-2t$, by the above claim, \textbf{Claim 4} holds.

For any two edges, say $e_1,e_2$, of $T_G$, let $e_1=u_1u_2$, $e_2=v_1v_2$.
Without loss of generality, we assume
$d_{T_G}(u_1,v_2)=\max\{ d_{T_G}(u_i,v_j): 1\leq i, j\leq 2\}$ where
$d_{T_G}(u,v)$ denote the distance between $u$ and $v$
in $T_G$.  As $T_G$ is a tree, the (unique) $u_1-v_2$ geodesic,
say $P$, in $T_G$ must contains edges $e_1,e_2$. Moreover,
it is easy to show $P$ is also an unique $u_1-v_2$ geodesic in $G$,
so $c(e_1)\neq c_(e_2)$ under any strong rainbow coloring.

By \textbf{Claim 3} and \textbf{Claim 4},
the equality holds if and only if $G\in \overline{\mathcal{G}}_t$.
Then our result holds.
\qed
\end{pf}

Next we give an application to Theorem \ref{thm2}, we consider the
strong rainbow connection numbers of line graphs of connected cubic
graphs. Recall that the $line~graph$ of a graph $G$ is the graph
$L(G)$ whose vertex set $V(L(G))=E(G)$ and two vertices $e_1$, $e_2$
of $L(G)$ are adjacent if and only if they are adjacent in $G$. The
star, denoted $S(v)$, at a vertex $v$ of graph $G$, is the set of
all edges incident to $v$. Let $\langle S(v)\rangle$ be the subgraph
of $L(G)$ induced by $S(v)$, clearly, it is a clique of $L(G)$. A
$clique~decomposition$ of $G$ is a collection $\mathscr{C}$ of
cliques such that each edge of $G$ occurs in exactly one clique in
$\mathscr{C}$.  An $inner~vertex$ of a graph is a vertex with degree
at least two. For a graph $G$, we use $\overline{V_2}$ to denote the
set of all inner vertices of $G$. Let $\mathscr{K}_0=\{\langle
S(v)\rangle: v\in V(G) \}$, $\mathscr{K}=\{\langle S(v)\rangle: v\in
\overline{V_2} \}$. It is easy to show that $\mathscr{K}_0$ is a
clique decomposition of $L(G)$ and each vertex of the line graph
belongs to at most two elements of $\mathscr{K}_0$. We know that
each element $\langle S(v)\rangle$ of ${\mathscr{K}_0}\setminus
{\mathscr{K}}$, a single vertex of $L(G)$, is contained in the
clique induced by $u$ that is adjacent to $v$ in $G$. So
$\mathscr{K}$ is a clique decomposition of $L(G)$.

\begin{cor}\label{cor1} Let $L(G)$ be the line graph of a connected
cubic graph with $n$ vertices, then $src(L(G))\leq n$.
\end{cor}
\begin{pf} Since $G$ is a connected cubic graph, each vertex is an
inner vertex and the clique $\langle S(v)\rangle$ in $L(G)$
corresponding to each vertex $v$ is a triangle. We know that
$\mathscr{K}=\{\langle S(v)\rangle: v\in \overline{V_2} \}=\{\langle
S(v)\rangle: v\in V \}$ is a clique decomposition of $L(G)$. Let
$\mathcal {T}=\mathscr{K}$. Then $\mathcal {T}$ is a set of $n$
edge-disjoint triangles that cover all edges of $L(G)$. As there are
$3n$ edges in $L(G)$, by Theorem \ref{thm2}, we have $src(L(G))\leq
3n-2n=n$.\qed
\end{pf}

\end{document}